\documentclass{article}
\usepackage{amsfonts,amssymb,amsmath}

\usepackage[cp1250]{inputenc}

\usepackage{graphicx}
\usepackage{amsmath, amsthm, amssymb}
\usepackage{amsfonts}
\usepackage{enumerate}
\usepackage{verbatim}
\usepackage{color}

\newcommand{\tluste}[1]{\mbox{\mathversion{bold}$ #1 $}}
\newcommand{\A}[0]{{\tluste{A}}}
\newcommand{\B}[0]{{\tluste{b}}}
\newcommand{\X}[0]{{\tluste{x}}}
\newcommand{\Y}[0]{{\tluste{y}}}
\newcommand{\R}[0]{{\mathbb{R}}}

\newcommand{\IR}[0]{{\mathbb{IR}}}

\newcommand{\ol}[1]{\mbox{$\overline{{#1}}$}} 
\newcommand{\ul}[1]{\mbox{$\underline{{#1}}$}} 
\newcommand{\sgn}{\mathop{\rm sgn}\nolimits}
\newcommand{\diag}{\mathop{\rm diag}\nolimits}

\newcommand{\www}{\mathop{\rm w}\nolimits}

\begin{document}
\begin{center}

{\LARGE\bf
Computing Enclosures of Overdetermined Interval Linear Systems
}\bigskip

{\large
Jaroslav Hor\'a\v cek$^{1}$ and Milan Hlad\'ik$^{1,2}$
}\medskip

$^1$ Charles University, Faculty of Mathematics and Physics, Department of Applied
Mathematics, Malostransk\'e n\'am. 25, 118 00, Prague, Czech Republic
{\tt
horacek@kam.mff.cuni.cz, hladik@kam.mff.cuni.cz
}
\newline

$^2$ University of Economics, Faculty of Informatics and Statistics,
n\'am. W. Churchilla 4, 13067, Prague, Czech Republic

\end{center}


\begin{abstract}
This work considers special types of interval linear systems - overdetermined systems. Simply said these systems have more equations than variables. The solution set of an interval linear system is a collection of all solutions of all instances of an interval system. By the instance we mean a point real system that emerges when we independently choose a real number from each interval coefficient of the interval system. Enclosing the solution set of these systems is in some ways more difficult than for square systems. The main goal of this work is to present various met\-hods for solving overdetermined interval linear systems. We would like to present them in an understandable way even for nonspecialists in a field of linear systems. The second goal is a numerical comparison of all the methods on random interval linear systems regarding widths of enclosures, computation times and other special properties of methods.    
\end{abstract}

$\\$

\noindent Research is supported 
by the Grant Agency of the Charles University (GAUK) grant no. 712912. \newline
Jaroslav Hor\'a\v cek was partially supported by the Czech Science Foundation under the contract
201/09/H057. \newline
M.\ Hlad\'{\i}k was partially supported by the grant GA\v{C}R P403/12/1947.

$\\$

\section{Introduction}
Real-life problems can be described by different means -- by linear and nonlinear systems, by systems of difference and differential equations etc. Nevertheless the description can often be transformed to another one using only linear equalities (or inequalities). To take into account rounding errors or imprecisions of measurement or data we can use means of interval analysis. That is why interval linear systems are still in the focus of researchers. There are plenty of methods for enclosing the solution set of square interval linear systems -- systems in the form $\A x=\B$, where $\A$ is a square matrix (e.g. \cite{moore:introduction}, \cite{Neumaier1990}, \cite{NinKea1997}, \cite{Rump2010}). That is because square matrices can posses some advantageous properties. They can be diagonally dominant, positive definite, M-matrices and many more. And we know that many algorithms behave well for those cases. However sometimes we have to work with systems called \emph{overdetermined}, simply said, they consist of more equations than variables. That is why our favorite methods for square systems usually can not be applied. 

Fortunately there exist some methods for solving those systems -- Gaussian elimination, classical iterative methods, Rohn method, least squares method or linear programming. 
The main goal of this work is to present the overview of existing methods for computing enclosures of solution sets of overdetermined interval linear systems. To the best of our knowledge, that has not been done for overdetermined systems yet. We would like to explain how these methods work in a brief but understandable way even for researchers from various fields not so familiar with interval linear systems. 
We also would like to mention some pitfalls and specialties connected with these methods. Some of them behave in an useful way -- they are able to compute a very narrow enclosure or an interval hull or they can reveal unsolvability of an interval overdetermined system. On the other hand, many of them can not decide whether a solution of an interval system is unbounded or whether the system is unsolvable. Moreover sometimes if a system has certain properties, they are not able to return any meaningful result.
It is interesting to observe how some efficient methods fail when the radii of intervals change and, on the other hand, to see how some simple, one could say "stupid" methods rule. After introducing the methods we provide the result of exhaustive numerical comparison of these methods concerning the speed, enclosure tightness and some other properties.
As always it holds that not every method is useful for every problem case. 
When describing the methods we would like to point out the cases when they are useful.

\section{Basic notation and definitions}
It is favourable to start with the basic notation. Here we provide only a small part of it which we will use subsequently. For deeper introduction to interval arithmetics do not hesitate to see \cite{moore:introduction}. 
We denote interval structures in boldface ($\A, \B$). Point real structures are denoted in regular type ($A, b$). Here we work only with closed real intervals $[c,d],$ where $c \leq d$. The set $\IR$ stands for the set of all real closed intervals. An interval $\X$ can be defined in two ways. The first one is by upper and lower bound $\X=[\ul{x}, \ol{x}]$. The next one is by midpoint and radius $\X= \langle x_c, x_\Delta \rangle$. With intervals we can build more complex structures -- interval matrices (of which vectors are special case).  Their definition is analogous to the definition of a real interval ($\A=[\ul{A}, \ol{A}]$ or $\A = \langle A_c, A_{\Delta}\rangle$). The matrix $A_c$ is called the \emph{midpoint matrix} and $\A_\Delta$ is called the \emph{radius matrix}. In the section about comparison of methods we will use the notion \emph{width} of an interval $\X = [\ul{x}, \ol{x}]$ defined as $\www(x) = \ol{x} - \ul{x}$. 
For every vector $x \in \R^n$ we define its sign vector 
$\sgn(x) \in \{\pm 1\}^n$ as

\begin{equation}
(\sgn x)_i = \left\{ \begin{array}{r l}
1, & \textrm{ if $x_i \geq 0 $,}\\
-1, & \textrm{ if $x_i < 0 $. }\\
\end{array} \right.\nonumber
\end{equation}  

\noindent For a given vector $x \in \R^{n} $ we denote

\begin{displaymath}
D_x = \diag( x_1, \dots, x_n) = 
\left( \begin{array}{cccc}
x_1 & 0 & \ldots & 0 \\
0 & x_2 & \ldots & 0 \\
\vdots & \vdots & \ddots & \vdots \\
0 & 0 & \ldots & x_n
\end{array} \right).
\end{displaymath}

\noindent We continue with a definition of an overdetermined interval linear system.  

\newtheorem*{intlinSystem}{Definition (Overdetermined interval linear system)}
\begin{intlinSystem}
Let us have an interval matrix $\A \in \IR^{m \times n}$, where $m >n$  and an interval vector $\B \in \IR^{m}$. We call
$$\A x = \B $$
an overdetermined interval linear system (OILS). 
\end{intlinSystem}
\noindent Simply said an overdetermined interval linear system is just a general interval linear system that consist of more equations than variables. 

When we talk about the interval linear systems we have to mention what do we mean by the solution of an interval linear system (ILS).

\begin{figure}[h]
		\centering
			\includegraphics[width=6cm]{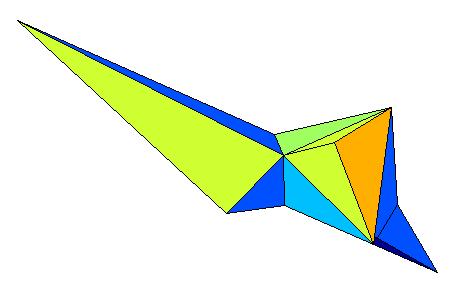}
		\caption{Graphical view of a solution of an interval linear system with 3 variables}
		\label{fig:solset}
\end{figure}

\newtheorem*{ilsSolution}{Definition (Solution set of ILS)}
\begin{ilsSolution}
The solution set $\Sigma$ of an interval linear system $\A x = \B$ is 
	$$\Sigma = \{ \, x \ | \ Ax=b \ \textrm{for some} \ A \in \A, b \in \B \, \}.$$	
\end{ilsSolution}
\noindent In another words it is a collection of all solutions of all instances of an interval linear system. By an instance we mean a point-real system that we get when we independently choose a real number from each interval coefficient of an interval system. We have to mention that this approach is different from the least squares approach. If no instance of the system has solution, we call the whole interval system \emph{unsolvable}. To give a reader better conception of an interval solution set, 
on the picture~\ref{fig:solset} we show the solution set of an interval system

$\\$
\begin{displaymath}
\begin{array}{ccccccc}
[-5,10] \, x & + & 10\, y & + & [15,20]\, z & = & [50,100],\\
10 \, x & + & -5\, y & + & [5,15]\, z & = & [-50,50], \\
10\, x & + & [10, 25]\, y & + & [-10,-5]\, z & = & [50,100].\\
\end{array}
\end{displaymath}
$\\$

\noindent The solution set is generally a polyhedral set, not necessarily convex. Never\-theless it is convex in each orthant. As we can see, this set is quite difficult to be described. There exist a lot of possibilities of description. One of them  is to find an n-dimensional box (aligned with axes) as tight as possible that contains the solution set. We call this box the \emph{interval hull}. Unfortunately it can be proved that computing the interval hull is NP-hard. That is why we look for a box a little bit wider, still as narrow as possible, that contains the interval hull. We call this box the \emph{interval enclosure}. There is a large variety of methods we can use to compute an interval enclosure. We have already built sufficient terminology, therefore we are able to move on to description of methods.

\section{Gaussian elimination (GE)}
The first thing that can come to our mind is the Gaussian elimination. The interval GE was proposed by Hansen in \cite{hansen:ols}. The idea is pretty the same as for point real systems. We present a slightly modified version.  We eliminate rows using interval operations with the only difference that we eliminate the system $(\A | \B)$ to the shape 

\begin{displaymath}
(\A \ | \ \B) \sim \ldots \sim 
\left(
\begin{array}{cc|c}
 \tluste{C} & \tluste{d}& \tluste{e}\\
 \tluste{0} & \tluste{f}& \tluste{g}\\
\end{array}
\right),
\end{displaymath}
where $\tluste{C}$ is an $(n-1) \times (n-1)$ interval matrix in the Row echelon form (REF) with $[1,1]$ intervals on pivot positions, $\tluste{d}, \tluste{e}$ are $(n-1) \times 1$ interval vectors, $\tluste{0}$ is an $(m-n+1) \times (m-n+1)$ matrix composed of the intervals $[0,0]$ and $\tluste{f}, \tluste{g}$ are $(m-n+1)\times 1$ interval vectors.   
\newline

\noindent One thing is often not clear to people familiar with interval arithmetics when Gaussian elimination is explained to them. And that is why on pivot positions there are $[1,1]$ intervals and why there are $[0,0]$ intervals in $\tluste{0}$ matrix. This might not be completely obvious. However the reason is quite simple -- we can think of working with pivots and elements beneath as performing the elimina\-ting operations on them in all instances of an interval system separately (after elimination pivots in instances are equal to 1, elements under the pivots are equal to 0). When we assemble the elimination result, we know that we get $[1,1], [0,0]$ intervals on these positions in all these instances, therefore there is no need to overestimate them using interval arithmetics. For the other coefficients we have to use interval arithmetics to gain their interval enclosures.
\newline

\noindent Now we can realize that vectors $\tluste{f}, \tluste{g}$ form $m-n+1$ interval equations in the shape $$\tluste{f}_i x_n = \tluste{g}_i \quad \textrm{for} \ i=1, \ldots, (m-n+1).$$ The solution of these equations is intersection of all intervals $\tluste{g}_i / \tluste{f}_i$ -- we get an enclosure of the variable $x_n$. If the intersection is empty, then the system has no solution. Nonetheless, if the intersection is unbounded, it can either mean that the solution set of the system is unbounded or large overestimation due to large number of interval operations occurred. 
The enclosures for the other variables can be obtained using backward substitution as we know it (only computing in interval arithmetics).

This algorithm works only for very small $m \times n$ systems $(n \sim 4)$. For larger systems we get a huge overestimation as we can see in table \ref{tab:gegrow}.

\begin{table}[h]
\centering
\begin{tabular}{c|c|c}
variable     & midpoint & radius  \\
\hline
$x_1$ &  0.1479 & 227.6698 \\
$x_2$ &  15.2091 & 172.5929 \\
$x_3$ &  11.1031 & 68.4653  \\
$x_4$ &  9.7809 & 64.1056  \\
$x_5$ &  -8.8168 & 27.2234  \\
$x_6$ &  25.8164 & 25.8398  \\
$x_7$ &  -19.0444 & 30.4596\\
$x_8$ &  -22.0799 & 11.0313  \\
$x_9$ &  1.9649 & 12.1172  \\
$x_{10}$ &  -19.1817 & 11.6841 \\ 
$x_{11}$ &  -20.9670 & 1.9153 \\
$x_{12}$ & -4.6988 & 3.5407 \\
$x_{13}$ & 3.1223 & 4.3894 \\
\end{tabular}
\caption[GE -- grow]{The overestimation of the enclosure (random system $15\times13$) with random radii $< 10^{-3}$ caused by GE and backward substitution}
\label{tab:gegrow}
\end{table} 

For larger systems it is necessary to use preconditioning. We use the preconditioner mentioned by Hansen in \cite{hansen:ols}.
\begin{displaymath}
C = 
\left[ \begin{array}{cc}
A^{c}_1 & 0 \\
A^{c}_2 & I \\
\end{array} \right]^{-1},
\end{displaymath}
where $A^{c}_1$ consists of first $n$ rows of $A_c$ (approximation of a midpoint matrix of $\A$) and $A^{c}_2$ consists of the remaining $m-n$ rows of $A_c$. Operation $^{-1}$ is an approximate inverse operation. Of course a preconditioning usually causes an overestimation of resulting enclosure, but it often causes smaller enlargement compared to performing GE without preconditioning. We must also mention that preconditioning enlarges the solution set of the system, that is why unsolvable systems become solvable after preconditioning. Therefore if we want to check unsolvability we can not use preconditioning. If we use the elimination to Reduced row echelon form (RREF) and allow multiple right-hand side, we are able to compute an interval matrix containing $\A^{-1}$. For a proper definition of the interval inverse matrix see \cite{rohn:handbook}.

\section{Iterative methods}
There are many interval modifications of methods for square non-interval systems -- Jacobi, Gauss-Seidel, Krawczyk, Hull method. But we usually can not use them when dealing with overdetermined systems. Moreover without preconditioning these methods often work really badly. However we have one possibility, if we want to use our known iterative methods we can use the preconditioning mentioned in the section about Gaussian elimination. We can see that after this preconditioning the new system 
is usually of the shape
\begin{displaymath} 
\left( \begin{array}{c}
\tluste{I^{\sim}} \\
\tluste{0^{\sim}}
\end{array} \right), 
\end{displaymath} 
where $\tluste{I^{\sim}}$ is an $n \times n$ interval matrix "wrapping" the point eye matrix and $\tluste{0^{\sim}}$ is an $(m-n)\times n$ interval matrix with small intervals containing zero. Then we can use our favorite iterative method for the upper square $n \times n$ subsystem. Let us choose Jacobi method for example. We start with some initial enclosure $\X^{(0)}$ and iteratively "sharpen" this enclosure using the following formula (in each step we simply express the variable $x_i$ from $i$-th equation, for more information see \cite{moore:introduction}) 

\begin{equation}
\label{jacobivzorec}
\X^{*}_i = \frac{1}{\tluste{A}_{ii}}( \B_i - \sum_{j \neq i} \tluste{A}_{ij} \X^{(k)}_j) \quad \textrm{for}  \ i= (1, \ldots, n), \ \textrm{in} \ (k+1)\textrm{-th step}.
\end{equation}

\noindent After every iteration we provide intersection with the old enclosure to obtain
$$ \X^{(k+1)}_i = \X^{(k)}_i \cap \X^{*}_i \quad \textrm{for}  \ i= (1, \ldots, n). $$
Gauss-Seidel method differs in only one thing -- we do not wait and compute intersection after every computation of $\X^*_i$ and immediately use the new value for further computing. This leads to less iteration steps.
\newline

\noindent The problem is that when using these methods, intervals containing zero are are not allowed on a diagonal (since there is the division with $\A_{ii}$). That case happens if the radii of the original matrix are "large" (even $10^{-2}$).
This method looks rather simple. We loose some information due to precoditioning and chopping the last $m-n$ rows of the original system. But we can realize that e.g. Jacobi method can be effectively parallelized. And if there is effective way to determine $\X^{(0)}$, then this method can be used for very fast sharpening of $\X^{(0)}$, for example in Constraint Programming as a sharpening step in between the run of another sharpening methods.  

\section{Rohn method}
We would like to mention the method introduced by Rohn. Due to lack of space we will introduce it only briefly. For more information and theoretical insight one can take a look in  \cite{rohn:enclosing}. The base of the method lies on the following theorem.

\newtheorem*{rohn}{Theorem (Rohn)}
\begin{rohn}
Let $\A x = \B$ be an IOLS with a solution set $\Sigma$ ($\A$ is an $m \times n$ matrix).
Let $R$ be an arbitrary real $n \times m$ matrix, let $x_{0}$ and $d > 0$ be arbitrary $n$-dimensional real vectors such that
$$ Gd + g < d, $$
\noindent where
$$ G = |I - R A_c| + |R|A_\Delta $$
\noindent and
$$ g = | R(A_c x_0 - b_c) | +  |R|(A_\Delta|x_0| + b_\Delta). $$
\noindent Then 
$$ \Sigma \subseteq [x_0 - d, x_0 + d]. $$
\end{rohn}
$\\$
\noindent The question is how to find the vector $d$, the matrix $R$ and the vector $x_o$. To compute $d$, we can, for example, rewrite the inequality as 
$$ d = Gd + g + f, $$
for some small vector $f >0$. Then start with $d=0$ and iteratively refine $d$. This algorithm will stop after a finite number of steps if the spectral radius of $G$ is less than 1. Otherwise we do not know what happens. During practical testing with random systems with radii of intervals close to $0.1$ the vector $d > 0$ was rarely found.  
\newline

\noindent We still have to determine $x_0$ and $R$. For the start we can take $$ x_0 \approx R b_c,$$ $$ R \approx (A^T_c A_c)^{-1} A^T_c $$ but not necessarily. The theorem provides a clever instrument for iterative improvement of solution. We do not have to use only $A_c$ to compute $R$, we can take any (e.g. random) $A \in \A$, compute an enclosure and then intersect it with the old one.  We can repeat this process as many times we want and provide an iterative improvement of the enclosure. For smaller systems the iterative improvement works pretty well. The narrowing of initial enclosures is showed in the table \ref{tab:rohniteracnid}. It shows the ratios of enclosure widths returned by non-iterative Rohn method and iterative versions for 10, 100, 1000 iterations. To see how the ratios are computed, it is possible to take a look in subsection 8.1. In the further text we will call the first method \emph{basic Rohn method} and use the notion \emph{basic enclosure}. When we use the iterative improvement we will talk about the \emph{iterative Rohn method} and \emph{iterative enclosure}.

\begin{table}[h]
\centering
\begin{tabular}{c|c|c|c}
 system     & 10 it.& 100 it.& 1000 it.\\
\hline
$5 \times 3$     & 0.73   & 0.57   & 0.50  \\
$ 15 \times 10$  & 0.89  &  0.82  & 0.76 \\
$ 25 \times 21 $ & 0.94   &  0.90  & 0.87 \\
$ 35 \times 23 $ & 0.95   &  0.92  & 0.90 \\
$ 50 \times 35 $ & 0.97   & 0.94   & 0.93 \\
$ 70 \times 55 $ & 0.98   & 0.96   & 0.95 \\
$ 100 \times 87 $ & 0.98  & 0.97   & 0.97 \\ 
\end{tabular}
\caption[Rohn -- ratios initial enclosure]{Rohn method -- ratios of widths of basic enclosures and iterative enclosures (10, 100 a 1000 iterations)}
\label{tab:rohniteracnid}
\end{table} 

\section{Least squares method}

We can also make use of the least squares formula $A^T A x = A^T b$. But when we write it in an interval way 
$$\A^T \A x = \A^T \B$$
that indeed does not work because of two interval matrix multiplication. Even if we use some preconditioner $C$ and write
$$(C \A)^T (C\A) x = (C \A)^T \B$$
that does not work either.
Anyway we can use an equivalent expression of the least squares formula for point real systems. If we use it for interval systems we get
\begin{displaymath}
\left( \begin{array}{cc}
I & \A \\
\A^T & 0 \\
\end{array} \right)
\left( \begin{array}{c}
y\\
x\\ 
\end{array} \right)
=
\left( \begin{array}{c}
\B\\
0\\ 
\end{array} \right).
\end{displaymath} 
\newline
Now we have a square system and we can apply some suitable method for square systems.
We get the vector $(\Y, \X)^T \in \IR^{m+n}$ as the result and take $\X$ as the solution enclosure. It can be seen that the returned interval vector contains the solution of the interval least squares, therefore this method returns a solution even if the system is unsolvable. Another drawback is that if the original system is of size $m \times n$ we have to solve a new one of size $(m+n) \times (m+n)$ On the other hand this method computes very sharp interval enclosures. It is used in the method \verb|verifylss()| in INTLAB 6 \cite{rump:intlab} in combination with the interval Krawczyk method (described e.g. in \cite{moore:introduction}) for square systems.  
When solving a systems in the shown way, we can notice the dependencies in the new system. We can possibly use methods dealing with dependencies in interval linear systems (e.g. \cite{Hladik2012}, \cite{Popova2000}, \cite{Rump2010}).

\section{Linear programming (LP)}
\newtheorem*{oettliPrager}{Theorem (Oettli-Prager)}

It is possible to compute the interval hull by using linear programming. We can use the famous theorem by Oettli-Prager. It can be found with proof for example in \cite{rohn:linsys}. Due to lack of space, here we omit the proof.

\begin{oettliPrager}
Let us have an interval linear system $\A x = \B$.
\noindent Vector $x \in \R^{n} $ is a solution of this system ($x \in \Sigma$) if and only if
$$ | A_c x - b_c | \leq A_\Delta |x| + b_\Delta. $$
\end{oettliPrager}
\noindent Unfortunately we are still not able to use LP because of the absolute values. Now we demonstrate how to rewrite this problem using linear inequalities only. We can get rid of the first one by decomposing it in two cases 
\begin{eqnarray} 
\label{firstlp} 
A_c x - b_c  \ \leq  A_\Delta |x| + b_\Delta, \\
\label{firstlp2}
-( A_c x - b_c ) \ \leq  A_\Delta |x| + b_\Delta. 
\end{eqnarray}
The second absolute value can be rewritten with the use of knowledge of the orthant we currently "are" inside. The following holds 
$$|x| = D_z x, \ \textrm {where} \ z=\sgn(x). $$
That gives a rise to the condition
\begin{equation}
0 \leq D_z x. 
\label{secondlp}
\end{equation}

\noindent Now we can see, that for every orthant conditions (\ref{firstlp}), (\ref{firstlp2}) and (\ref{secondlp}) form a system of linear inequalities. Therefore we can use linear programming. Unfortunately we have to use $(2^n \times 2n)$ linear programming problems (in each coordinate we compute the upper and lower bound). That is obviously too much computing. 
However we can compute the enclosure of the solution of the system with some suitable method (least squares, Rohn) and than apply linear programming to only those orthants where this enclosure lies. This approach is often much more faster. The table \ref{tab:linprogtime} illustrates a large speedup when solving the system with different method before applying LP. The sign '-' means that we omitted the testing because of enormous computing time when compared to  LP with presolving.

\begin{table}[ht]
\centering
\begin{tabular}{c|c|c}
system & time LP & time LP presolved\\
\hline
$5 \times 3$  & 6 sec & 1 sec \\
$ 9 \times 5$ &  43 sec  & 1.68 sec\\
$ 13 \times 7$ &  5 min & 3.59 sec\\
$ 15 \times 9$  &  28 min & 4.1 sec \\
$ 25 \times 21$ & - & 13 sec\\
$ 35 \times 23$  & - & 19 sec \\
$ 45 \times 31$  & - & 43 sec \\
$ 55 \times 35$  & - & 1 min  \\
$ 73 \times 55$  & - & 9 min  \\
\end{tabular}
\caption[LP times some]{Comparison of times of LP in all orthants and LP with presolving with a different method (here verifylss)}
\label{tab:linprogtime}
\end{table}
  
\section{Comparison of methods}
The another goal of this paper is a numerical comparison of methods for enclosing solutions of overdetermined systems. As far as we are concerned for these methods it  has not been done yet. All subsequent tests were computed using the following hardware and software:
\begin{itemize}
\item Processor -- AMD Phenom(tm) II X6 1090T 
\item Memory -- 15579 MB
\item Matlab R2010b
\item Intlab 6 (see \cite{rump:intlab})
\item Versoft 10 (see \cite{versoft}) for linear programming
\end{itemize}
It is difficult to imagine rectangular matrices in some special shape, therefore we test the methods on random matrices. In this work we use four parameters:
\begin{itemize}
\item \emph{maximum radius}
\item \emph{midpoint range}
\item \emph{stopping parameter}
\item \emph{maximum number of iterations}
\end{itemize}

Most of these parameters are used for generating random interval systems with certain properties. The key parameter is \emph{maximum radius} of interval coefficients of an interval linear system. When this parameter is for example $-4$, that means that for every interval coefficient $\Y = \langle y_c, y_{\Delta} \rangle$ in the system
$y_{\Delta} \leq 10^{-4}$ holds. It is usually negative, because radii greater than 1 will lead to singularities. The radii in a system are not the same, they are chosen randomly with respect to the maximum radius parameter. Another important parameter is the \emph{midpoint range}, which specifies the range of midpoints of intervals. The methods may have different properties when working with all intervals that have relatively close midpoints and between the case when the midpoints are quite far like in vector $ (10, 1, 11234, 0.1)^T$. We will test on two cases -- the first one is with midpoints uniformly chosen from $[-25, 25]$ and the second one uniformly chosen from $[-1000, 1000]$. Some methods need small real positive number $\epsilon$ as a \emph{stopping parameter} (some methods need a small positive vector, but we can use $f=(\epsilon, \epsilon, \ldots, \epsilon)^T$). We will choose $\epsilon = 10^{\,( \textrm{maximum radius} - 2)}$. After performing some practical tests this seems to be a reasonable choice. Some iterative methods use parameter \emph{maximum of iterations}. It is a safety precaution for the case the stopping parameter did not work. We use the value 20 for all the tested cases. We will test systems up to size about $m=1000$ which equals to a thousand of linear equations. We chose this upper bound on size of tested systems after consulting with a colleague from the Faculty of Civil Engineering. These are maybe the largest linear systems they need to solve for many technical purposes. If a solution of larger system is needed, then it is usually better to split the  problem into parts, solve the parts separately and then assemble the solution.  
For the sake of clarity we establish the following identifiers of methods:

\begin{itemize}
\item \emph{GSpre} - Gauss-Seidel iterative method with preconditioning
\item \emph{GEpre} - Gaussian elimination with preconditioning
\item \emph{Rohn} - basic Rohn method
\item \emph{Verifylss} - method verifylss from toolbox Intlab by Rump
\item \emph{LPver} - linear programming from toolbox Versoft by Rohn
\end{itemize}

\subsection{Enclosure width comparison}
We are going to compare the widths of enclosures of all the mentioned methods according to test enclosures returned by a different method.
For an $n$-dimensional enclosure vector $\X=(\X_1, \X_2, \ldots, \X_n)^T$ and 
$\X^{test}$ as a test enclosure we compute the ratio
$$ \mathop{\rm ratio}(\X, \X^{test}) = \frac{1}{n} \sum^{n}_{i=1} \frac{ \www(\X^{}_i)}{ \www(\X^{test}_i) }. $$
\noindent Clearly, the lower the ratio, the better the enclosure. 

We divide the testing in two phases. For larger systems it is very time consuming to compute many enclosures using \emph{LPver}, moreover for small radii of intervals the method returns $NaN$ solution. That is why we tested the methods on small systems first and the enclosures were compared with \emph{LPver} results and in second phase they were compared to \emph{Verifylss} enclosures.  
First let us show the table \ref{tab:widthssmall3} where the ratios of small systems are displayed.

\begin{table}[h]
\centering
\begin{tabular}{c|c|c|c|c}
system & $GSpre$ & $GEpre$ & $Rohn$ &  $Verifylss$ \\
\hline 
$5 \times 3$   & 2.9968 & 2.9969 & 1.2347  & 1.1893 \\
$15 \times 9 $ & 7.7057 & 7.7069 & 1.1601 & 1.1500 \\
$ 35 \times 23 $ & 8.3591 & 8.3602 & 1.1276 & 1.1249 \\
$55 \times 35 $ & 11.5026 & 11.5064 & 1.1336 & 1.1331 \\
$73 \times 55$ & 17.6528 & 17.6624 & 1.0828 & 1.0848 
\end{tabular}
\caption[Small radii of systems widths]{Ratios of widths of enclosures, \emph{LPver} is 1 (maximum radii $\leq 10^{-3}$)} 
\label{tab:widthssmall3}
\end{table}

We can see that \emph{GEpre} and \emph{GSpre} are doing really badly. And that their enclosure ratio is of about the same value. We attribute that to the use of the same preconditioner. \emph{Rohn} and \emph{Verifyllss} provide good results. \emph{Verifylss} is a little bit better, but the difference is slowly disappearing with the grow of a system size. That is confirmed by the second table \ref{tab:widthsbig5}. We made the radii smaller ($\leq 10^{-5}$) since \emph{GEpre} often returned infinite solution due to a large number of interval operations.

\begin{table}[h]
\centering
\begin{tabular}{c|c|c|c}
system & $GSpre$ & $GEpre$ & $Rohn$ \\
\hline 
$100 \times 45$ & 21.1668 & 21.1668 & 1.0237 \\
$100 \times 87$ & 9.2420 & 9.2420 & 1.0062 \\
$180 \times 125$ & 26.6744 & 26.6744 & 1.0064 \\
$180 \times 170$ & 15.0262 & 15.0262 & 1.0020 \\
$290 \times 190$ & 37.6976 & 37.6977 & 1.0044 \\
$290 \times 260$ & 26.9361 & 26.9361 & 1.0018 \\
$380 \times 275$ & 40.5159 & 40.5160 & 1.0028 \\
$380 \times 360$ & 32.9500 & 32.9504 & 1.0009 \\
$500 \times 350$ & 66.9764 & 66.9771 & 1.0022 \\
$500 \times 470$ & 34.2512 & 34.2516 & 1.0007 
\end{tabular}
\caption[Small radii of systems widths]{Ratios of widths of enclosures, \emph{Verifylss} is 1 (maximum radii $\leq 10^{-5}$)} 
\label{tab:widthsbig5}
\end{table}

\subsection{Enclosure computation time comparison}
We now compare the computation times of methods - \emph{GSpre}, \emph{GEpre}, \emph{Rohn}, \emph{Veri\-fylss}. Midpoint matrices and vectors are chosen uniformly from the interval [-1000, 1000]. We provide two tables \ref{tab:alltimes3} and \ref{tab:alltimes5} for maximum radius $\leq 10^{-3}$ and $\leq 10^{-5}$
respectively. They display the average computation times of all methods for systems of different sizes. Time is measured with Matlab functions \emph{tic} and \emph{toc}. The results are displayed only to the fourth decimal number. More digits are beyond the precision of \emph{tic} and \emph{toc}.

\begin{table}[h]
\centering

\begin{tabular}{c|c|c|c|c}
system & $GSpre$ & $GEpre$ & $Rohn$ & $Verifylss$ \\
\hline 
$5 \times 3$ & 0.0113 & 0.0117 & 0.00065 & 0.0021 \\
$15 \times 13$ & 0.0367& 0.0685 & 0.00086 & 0.0023 \\
$35 \times 23$ & 0.0636 & 0.2840 & 0.0012 & 0.0036 \\
$50 \times 35$ & 0.0983 & 0.5806 & 0.0019 & 0.0052 \\
$100 \times 87$ & 0.3063 & 2.4631 & 0.0115 & 0.0169 \\
$200 \times 170$ & 0.6632 & 9.8570 & 0.0523 & 0.0761 \\
$380 \times 275$ & 1.3253 & 34.5861 & 0.1765 & 0.3257 \\
$500 \times 470$ & 2.9486 & 69.6341 & 0.6235 & 0.9154 \\ 
\end{tabular}
\caption[All times -3]{Average times (in seconds) of all methods for maximum radius $\leq 10^{-3}$} 
\label{tab:alltimes3}
\end{table}

\begin{table}[h]
\centering

\begin{tabular}{c|c|c|c|c}
system & $GSpre$ & $GEpre$ & $Rohn$ & $Verifylss$ \\
\hline 
$5 \times 3$ & 0.0113 & 0.0118  & 0.00066 & 0.0021 \\
$15 \times 13$ & 0.0356 & 0.0687 & 0.00086 & 0.0023 \\
$35 \times 23$ & 0.0608 & 0.2836 & 0.0012 & 0.0036 \\
$50 \times 35$ & 0.0911 & 0.5809 & 0.0019 & 0.0053 \\
$100 \times 87$ & 0.2416 & 2.5224 & 0.0117 & 0.0171\\
$200 \times 170$ & 0.5718 & 10.6842 & 0.0561 & 0.0760 \\
$380 \times 275$ & 1.2515 & 37.0210 & 0.1924 & 0.3427\\
$500 \times 470$ & 2.5774 & 75.5171 & 0.6616 & 0.9761 \\ 
\end{tabular}
\caption[All times -5]{Average times in seconds of all methods for maximum radius $\leq 10^{-5}$} 
\label{tab:alltimes5}
\end{table}

We can see that the methods do not show any remarkable sensibility to the changes of radii of intervals if we keep the \emph{stopping criterion} equal to (\emph{maximum radius} $- 2$). There was also no extraordinary change when lowering the \emph{midpoint range} to something much smaller (e.g. $[-25,25]$). The complete time winner is \emph{Rohn}, the second is \emph{Verifylss}. The method \emph{GEpre} is doing really badly. We shall look more carefully on the comparison between \emph{Rohn} and \emph{Verifylss}. The method \emph{Rohn} is almost always faster. According to the table \ref{tab:widthsbig5} we can object that it pays off to take a little bit more time and compute more exact enclosure using \emph{Verifylss}. But let us look closer on the following table \ref{tab:rohnvsverifylss}. It displays the ratios of computation times of \emph{Rohn} and \emph{Verifylss} for systems of broader spectrum of sizes. We did the tests for many random systems, here we selected some interesting ones. It is possible to see, that for $m \times n$ systems, where $m$ is much bigger than $n$,  \emph{Rohn} is much more faster. For example on system $278 \times 35$ we can call the iterative \emph{Rohn} with 13 iterations within the same time \emph{Verifylss} needs to do the computation. And we can possibly sharpen the enclosure. 
The question whether this procedure pays off is answered by the table \ref{tab:rohniteracnid}. For relatively small systems it definitely pays off. For larger systems we might get not so significant improvements and it is safer to use \emph{Verifylss}.    

\begin{table}[h]
\centering

\begin{tabular}{c|c}
system & $\frac{t(Verifylss)}{t(Rohn)}$ \\
\hline

$80 \times 41$   &  2.8 \\
$189 \times 166$ &  1.5 \\
$278 \times 35 $ &  13.8\\
$377 \times 319$ &  1.6\\
$525 \times 285$ &  2.7 \\
$712 \times 271$ &  4.5\\
$807 \times 68$  &  46.5\\
$894 \times 8 $  &  423.7\\
$894 \times 797$ &   1.5\\
$906 \times 128$ &   19.1\\
$978 \times 235$  &  9.4\\
$1000 \times 663$ &  2.0
\end{tabular}
\caption[Rohn versus Verifylss]{Ratios of \emph{Verifylss} and \emph{Rohn} times} 
\label{tab:rohnvsverifylss}
\end{table}

\section{Conclusion} 
We introduced many methods for solving overdetermined interval linear systems -- Gaussian elimination, Iterative methods, Rohn method, Least square method and Linear programming. We mentioned some advantageous and unfavourable properties of these methods. At the end of this work we showed the comparison of all mentioned methods. As usually it is hard to tell which method is the best. Linear programming is suitable when we want to compute the interval hull of an interval system and have enough time. When our system has really small coefficients or our system is in the shape of "noodle", it is advantageous to use Rohn method otherwise it is favourable to use the least square method for fast computing of a sharp interval enclosure. Iterative methods can be used when we already know some enclosure and want to make try to sharpen it (e.g. as a part of some constraint solver) and then passing it as an input to ano\-ther method e.g LP. Gaussian elimination seems to work really bad, but when we consider the solvability of the system, it might give us a valuable information about the nonexistence of solution. 
\newline

\noindent 
Computing enclosures of solutions of overdetermined interval linear systems can occur as a subproblem of many computational tasks (constraint satisfa\-ction problem). And also the information about solvability of a system is often important (system validation, technical computing). There is still much work to do in this area. We believe that even faster and sharper algorithms can be developed. Simultaneously we are working on methods checking the solvability of overdetermined interval linear systems. Hopefully sufficient or necessary conditions for solvability or unsolvability of OILS can be found. The preconditioning is used for many methods, it would be useful to test different types of preconditioners. Some methods use probabilities, e.g. iterative Rohn method. There is a question whether we can effectively derandomize these methods.

\bibliographystyle{abbrv}
\bibliography{literatura}

\end{document}